 \documentclass[11pt,reqno]{amsart}

\usepackage{amssymb}

\usepackage{amsmath}

\makeatother

\newtheorem{theorem}{Theorem}
\newtheorem{lemma}{Lemma}

\newcommand{\lla}{\langle\!\langle}
\newcommand{\rra}{\rangle\!\rangle}
\renewcommand{\lll}{|\!|\!|}
\newcommand{\LASEP}{L}
\newcommand{\ASEP}{\lll_{-1,\lambda}}

\title[A Note on the Diffusivity of Finite-Range AEP's on $\mathbb Z$ ]
{A Note on the Diffusivity of Finite-Range Asymmetric Exclusion Processes on $\mathbb Z$}

\author{Jeremy Quastel $^{1}$}
\address{ Departments of Mathematics and Statistics,\\ University of
Toronto}
\email{quastel@math.toronto.edu}

\author{Benedek Valk\' o $^{1,2}$}
\address{ Departments of Mathematics and Statistics,\\ University of
Toronto}
\email{valko@math.toronto.edu}

\begin{document}

\footnotetext[1]{Supported by the Natural Sciences and Engineering Research Council of Canada.}
\footnotetext[2]{Partially supported by the Hungarian Scientific Research Fund grant K60708}

\begin{abstract}The diffusivity $D(t)$ of  finite-range asymmetric exclusion processes on $\mathbb Z$ with non-zero drift is expected to be of order $t^{1/3}$.  Sepp\"{a}l\"ainen and Bal\'azs recently proved this conjecture for the nearest neighbor case. We extend their results to general finite range exclusion by proving that the Laplace transform of the diffusivity is of the conjectured order.  We also obtain a pointwise upper bound for $D(t)$  of the correct
order.

\end{abstract}


\maketitle

\section{Introduction}
\label{intro}

A finite-range exclusion process on the integer lattice $\mathbb Z$ is a system of 
continuous time, rate one random walks with finite-range jump law $p(\cdot)$, i.e. $p(z)\ge 0$, and $p(z)=0$
for $|z|>R$ for some $R<\infty$, $\sum_z p(z)=1$, interacting via {\it exclusion}:  Attempted jumps to occupied
sites are suppressed.  We consider asymmetric exclusion process (AEP) with
non-zero drift, \begin{equation}\sum_z zp(z) = b\neq 0.\label{drift}\end{equation}  
The state space of the process is $\{0,1\}^\mathbb Z$. Particle 
configurations are denoted by $\eta$, with $\eta_x\in\{0,1\}$ indicating the absence, or presence, of a 
particle at $x\in \mathbb Z$. The infinitesimal generator of the process is given by
\begin{equation}
L f(\eta)=
\sum_{x,z\in \mathbb Z} p(z) \eta_x (1-\eta_{x+z})(f(\eta^{x,x+z}) -f(\eta) ) 
\end{equation}
where $\eta^{x,y}$ denotes the configuration obtained from $\eta$ by interchanging the
occupation variables at $x$ and $y$. 

Bernoulli product measures $\pi_\rho$, $\rho\in [0,1]$, with $\pi_\rho(\eta_x=1) = \rho$ form a one-parameter family of invariant measures for the process. The process starting from $\pi_0$ and $\pi_1$ are trivial and so we
consider the stationary process obtained by starting with $\pi_\rho$ for some $\rho\in (0,1)$.  
Although the fixed time marginals of this stationary process is easy to understand (since the $\eta$'s are just independent Bernoulli($\rho$) random variables), there are still lots of open questions about the full (space-time) distribution. Information about this process (and about the appropriate scaling limit) would be very valuable to understand such elusive objects as the Stochastic Burgers and the Kardar-Parisi-Zhang equations (see \cite{QV} for a more detailed discussion). 

We consider the two-point function,
\begin{equation}
S(x,t) = E[ (\eta_x(t)-\rho)( \eta_0(0)-\rho)],
\end{equation}
where the expectation is with respect to the stationary process obtained by starting from one of the invariant measures
$\pi_\rho$. $S(x,t)$ satisfies the sum rules (see \cite{PS})
\begin{equation} \label{MV}
\sum_x S(x,t) = \rho(1-\rho) = \chi,\qquad \frac{1}{\chi}\sum_x x S(x,t) = (1-2\rho)bt.
\end{equation}
The diffusivity
$D(t)$ is defined as
\begin{equation}\label{dee}
D(t) = ({\chi t})^{-1} \sum_{x\in \mathbb Z} (x - 
(1-2\rho) b t )^2 S(x,t).
\end{equation}
Using scaling arguments one  conjectures \cite{Sp},
\begin{equation}
S(x,t) \simeq  t^{-2/3} \Phi( t^{-2/3}(x- (1-2\rho) bt))
\end{equation}
for some  scaling function $\Phi$, as $t\to \infty$.  A reduced conjecture is that 
\begin{equation}\label{conj}
D(t) \simeq C t^{1/3},
\end{equation}
as $t\to \infty$.
Note that this means that the process has a superdiffusive behavior, as the usual diffusive scaling would lead to $D(t)\to D$. It is known that the mean-zero jump law would lead to this case, see \cite{V}.

If $f(t)\simeq t^{\rho}$ as $t\to \infty$ then  as $\lambda\to 0$,
\begin{equation}\label{tauberian}
\int_0^\infty e^{-\lambda t} t f(t) dt\simeq \lambda^{-(2 + \rho)}.
\end{equation}
 If $f$ satisfies (\ref{tauberian}) then we will say that $f(t) \simeq t^{\rho}$ {\it in the weak (Tauberian) sense}. Without some extra regularity
 for $f$
 (for example, lack of oscillations as $t\to\infty$),
 such a statement will not imply a strong version of $f(t)\simeq t^{\rho}$ as $t\to \infty$.  However,
 it does capture the key scaling exponent.  The weak (Tauberian) version of the conjecture (\ref{conj}) is  $\int_0^\infty e^{-\lambda t} t D(t) dt\simeq \lambda^{-7/3}$.

The first non-trivial bound on $D(t)$ was given in \cite{LQSY} using the so called \emph{resolvent approach}: the authors proved that $D(t) \ge C t^{1/4}$ in a weak (Tauberian) sense. (They also proved the bound  $D(t)
\ge C(\log t)^{1/2}$ in $d=2$, which was later improved to $D(t)\simeq C(\log t)^{2/3}$ in \cite{Y}.)
This result shows that the stationary process is indeed superdiffusive, but does not provide the conjectured scaling exponent $1/3$. 

The identification of this exponent  was given in the breakthrough paper of Ferrari and Spohn \cite{FS}. They treated the case of the totally asymmetric
simple exclusion process (TASEP) where the jump law is $p(1)=1$, $p(z)=0$, $z\neq1$. The focus of \cite{FS} is not the diffusivity, their main result is a scaling limit for the fluctuation at time $t$ of a randomly growing discrete one dimensional interface $h_t(x)$ connected to the equilibrium process of TASEP.  This random interface (the so-called \emph{height function}) is basically the discrete integral of the function $\eta_x(t)$ in $x$. The scaling factor in their result is $t^{1/3}$ and the limiting distribution is connected to the Tracy-Widom distribution. 
The proof of Ferrari and Spohn is through a direct mapping between TASEP and a particular last passage percolation problem, using a combination of results from \cite{BDJ}, \cite{J}, \cite{BR}, \cite{PS}.

The diffusivity $D^{TASEP}(t)$ can be expressed using the variance 
 of the height function (see \cite{QV}):
\begin{equation}\label{d-var}
D^{TASEP}(t) = (4\chi t)^{-1} \sum_{x\in \mathbb Z} Var(h_t(x))-4\chi |x-(1-2\rho)t| .
\end{equation}
This identity, the results of \cite{FS} and some additional tightness bounds would  imply the existence of the limit $D^{TASEP}(t) t^{-1/3}$ and even the limiting constant can be computed (see \cite{FS} and \cite{QV} for details). Unfortunately, the needed estimates are still missing, but from \cite{FS} one can at least obtain  a lower bound of the right order:
\begin{equation}\label{lbtasep}
D^{TASEP}(t)>Ct^{1/3} 
 \end{equation}
 with a positive constant $C$ (see \cite{QV} for the proof).

In \cite{QV} the resolvent approach is used to prove the following comparison theorem. 
\begin{theorem}[QV, 2006]\label{thm_univ}
Let $D_1(t), D_2(t)$ be the diffusivities of two 
 finite range exclusion processes in $d=1$
with non-zero drift.   There exists $0<\beta, C<\infty$ such that 
\begin{eqnarray}\nonumber
C^{-1}\int_0^\infty e^{-\beta \lambda t} t D_1(t) dt
&\le& \int_0^\infty e^{-\lambda t} t D_2(t) dt 
%\\ \nonumber &\le & 
\le C\int_0^\infty e^{-\beta^{-1} \lambda t} t D_1(t) dt\\\label{equiv}
\end{eqnarray}
\end{theorem}
Combining Theorem \ref{thm_univ} with 
 with   (\ref{lbtasep}) it was shown in \cite{QV} that
 \begin{theorem}[QV, 2006]\label{thm_tauberlb}
For any finite range exclusion process in $d=1$
with non-zero drift,
$D(t)\ge Ct^{1/3}$ in the weak (Tauberian) sense.
\end{theorem}
\noindent  \cite{QV}  also  converts the Tauberian bound into pointwise bound  in the  nearest neighbor case to get $D(t)\ge Ct^{1/3}(\log t)^{-7/3}$.

Just a few months later  Bal\'azs and Sepp\"al\"ainen  \cite{BS}, building on ideas of Ferrari and Fontes \cite{FF}, Ferrari, Kipnis and
Saada \cite{FKS}, and Cator and Groeneboom \cite{CG}, proved the following theorem:
\begin{theorem}[Bal\'azs--Sepp\"al\"ainen, 2006]\label{thm_BS}
For any  nearest neighbor asymmetric exclusion process in $d=1$ there exists a finite constant $C$  such that for all $t\ge 1$,
\begin{equation}\label{eq_BS}
C^{-1} t^{1/3}\le D(t) \le C t^{1/3}.
\end{equation}
\end{theorem}
Their proof uses refined and ingenious couplings to give bounds on the tail-probabilities of the distribution of the second class particle.

The aim of the present short note is to show how one can extend the results of \cite{QV} using Theorem  \ref{thm_BS}. Once we have the correct upper and lower bounds for $D(t)$ from (\ref{eq_BS}) in the nearest neighbor case, we can strengthen the statement of Theorem \ref{thm_tauberlb} using the comparison theorem:
\begin{theorem} \label{thm_tauber2}For any finite range exclusion process in $d=1$
with non-zero drift,
$D(t)=\mathcal{O}(t^{1/3})$ in the weak (Tauberian) sense:
there exists a constant $0<C<\infty$ such that
\begin{equation}\label{lu_bound}
C^{-1}\lambda^{-7/3}\le\int_0^\infty e^{-\lambda t} t D(t) dt \le C \lambda^{-7/3}.
\end{equation}
\end{theorem}
Getting strict estimates for a function using the asymptotic behavior of its Laplace transform usually requires some regularity and unfortunately very little is known qualitatively about $D(t)$. However, in our case (as noted in \cite{QV}), one can get an upper bound for $D(t)$ using an inequality involving $H_{-1}$ norms:
%
%
% by applying the following estimate (\cite{LY}) which holds for large enough $t$:
%\begin{equation}\label{eq_upper}
%t^{-1}\sum_x E[\int_0^t w (s) ds \int_0^t \tau_x w (s) ds] \le 4 \lll w\lll_{-1, t^{-1}}^2.
%\end{equation}
%Using (\ref{eq_upper}) with the Green-Kubo formula (\ref{GK}) we will get the following theorem.
%
\begin{theorem}\label{thm_upper}For any finite range exclusion process in $d=1$
with non-zero drift,
there exists $C>0$ such that for all $t\ge 1$,
\begin{equation}\label{dif_l_bound}
 D(t) \le C t^{1/3}.
\end{equation}
\end{theorem}
Note that in all these statements the constant $0<C<\infty$ depends on the jump
law $p(\cdot)$ as well as the density $0<\rho<1$. 

\section{Proofs}
\label{proofs}
\setcounter{equation}{0}

Theorem \ref{thm_tauber2} immediately follows from  Theorem \ref{thm_univ} using TASEP and a general finite range exclusion, together with Theorem \ref{thm_BS}. To prove Theorem \ref{thm_upper} one needs the Green-Kubo formula which relates 
the diffusivity to the time integral
of current-current correlation functions:
\begin{equation}
D(t) = \sum_z z^2 p(z)  +2\chi t^{-1} \int_0^t\int_0^s \lla w_0, e^{uL} w_0\rra duds. \label{GK}
\end{equation}
Here $w_x$ is the normalized microscopic flux
\begin{equation}
w_x=\frac1{\rho(1-\rho)} \sum_z p(z) (\eta_{x+z}-\rho)(\eta_x-\rho), 
\end{equation}
$L$ is the generator of the exclusion process and the inner product $\lla \cdot, \cdot \rra$ is defined for mean zero local functions $\phi, \psi$ as
\begin{equation}
\lla \phi, \psi\rra =E\left[ \phi  \sum_x\tau_x \psi\right],
\end{equation}
with $\tau_x$ being the appropriate shift operator.

(\ref{GK}) is proved in \cite{LOY} (in the special case $p(1)=1$,
but the proof for general AEP is the same.) 
A useful variant is obtained by taking the Laplace transform,
\begin{equation}\label{GK1}
\int_0^\infty e^{-\lambda t} tD(t) dt = \lambda^{-2} \left(
\sum_z z^2 p(z)+ 2\chi \lll w\ASEP^2 \right)
\end{equation}
where the $H_{-1}$ norm corresponding to $\LASEP$ is defined on a core of  local
functions by
\begin{equation}
\lll\phi\ASEP = \lla \phi, (\lambda - \LASEP)^{-1} \phi \rra^{1/2} .
\end{equation}
We also need the following inequality which has appeared in several similar versions in the literature (e.g.~\cite{LY},\cite{KL}).
\begin{lemma} \label{lem_upper} Let $w$ be the current for a finite range exclusion process in $d=1$
with non-zero drift. Then, 
\begin{equation}\label{eq_upper_}
t^{-1}\sum_x E[\int_0^t w_0 (s) ds \int_0^t w_x (s) ds] \le 12  \lll w_0\lll_{-1, t^{-1}}^2.
\end{equation}
\end{lemma}
\begin{proof} We will use the notation
\begin{eqnarray}
\lll \phi \lll^2=\lla \phi, \phi \rra,
\end{eqnarray}
note that the right hand side of (\ref{eq_upper_}) is equal to $t^{-1}  \lll \int_0^t w_0(s) ds \lll^2$.
Let $\lambda>0$ and $u_\lambda=(\lambda-L)^{-1} w_0$. Using Ito's formula together with the identity 
$
L u_{\lambda}=\lambda u_\lambda-w_0
$
we get
\begin{eqnarray}\label{mart}
u_{\lambda} (t)=u_{\lambda}(0)-\int_{0}^t (\lambda u_{\lambda}-w_0)(s) ds+M_t
\end{eqnarray}
where $M_t$ is a mean zero martingale with
\begin{eqnarray}
\lll M_t \lll^2&=&\int_0^t \lla u_{\lambda}(s), -L u_{\lambda}(s)   \rra ~ ds =t \lla u_{\lambda}, -L u_{\lambda}   \rra \nonumber\\
&=&  t \lla  u_\lambda, (\lambda-L) u_{\lambda}   \rra- t \lambda \lll u_{\lambda}\lll^2\label{quad}.
%\\
%&=&t \lll  w_0\lll^2_{-1,\lambda}   -t \lambda \lll u_{\lambda}\lll^2
\end{eqnarray}
Rearranging  (\ref{mart}), applying Schwarz Lemma and using stationarity 
\begin{eqnarray}
\lll \int_0^t w_0(s) ds\lll^2 &\le&8 \lll u_{\lambda}\lll^2+4 \lll M_t\lll^2+4 \lambda^2 \lll \int_0^t u_{\lambda}(s) ds\lll^2\label{l1}
\end{eqnarray}
To bound the last term of (\ref{l1}) we again use Schwartz Lemma and  stationarity to get 
\begin{equation}
\lll \int_0^t u_{\lambda}(s) ds\lll^2\le t\int_0^t \lll u_{\lambda}(s)  \lll^2 ds=t^2 \lll u_{\lambda} \lll^2 .\label{l2}
\end{equation}
Putting our estimates together, 
\begin{eqnarray}
 \lll \int_0^t w_0(s) ds\lll^2 &\le& (8-4 t\lambda+4 t^2 \lambda^2) \lll u_\lambda \lll^2+4 t \lla  u_\lambda, (\lambda-L) u_{\lambda}   \rra
\end{eqnarray}
%\\
%&\le &(8+4\lambda^2 t-4 \lambda t) \lll u_{\lambda} \lll^2+4 t \lll w_0 \lll^2_{-1,\lambda}
%\end{eqnarray}
Setting $\lambda=t^{-1}$ and dividing the previous inequality by $t$:
\begin{eqnarray}
t^{-1} \lll \int_0^t w_0(s) ds\lll^2&\le& 8 \lambda \lll u_\lambda \lll^2+4  \lla  u_\lambda, ( \lambda-L) u_{\lambda}   \rra \nonumber\\
&\le& 12  \lla  u_\lambda, ( \lambda-L) u_{\lambda}   \rra =12   \lll w_0 \lll^2_{-1,\lambda^{-1}}
,\nonumber
%\\&\le&4 \lll w_0 \lll^2_{-1,t^{-1}}
\end{eqnarray}
which proves the lemma.
\end{proof}

\begin{proof}[ Proof of   Theorem \ref{thm_upper}]
Theorem \ref{thm_tauber2} and identity (\ref{GK1}) gives
\begin{equation}
 \lll w_0 \lll^2_{-1,\lambda}\le C \lambda^{-1/3}\label{end}
\end{equation}
if $\lambda>0$ is sufficiently small. 
To bound $D(t)$ we use the Green-Kubo formula (\ref{GK}) noting that the second term on the right is equal to $\chi t^{-1} \lll \int_0^t w_0(s) ds \lll^2$. 
Using  Lemma \ref{lem_upper} and then (\ref{end}) to estimate this term we get that
\begin{eqnarray}
D(t)\le \sum _z z^2 p(z)+12 \chi  \lll w_0 \lll^2_{-1,t^{-1}}\le  C t^{1/3}
\end{eqnarray}
for large enough $t$.
\end{proof}

%\bibliography{bibtex}

\begin{thebibliography}{McK76}

%CMJ style

%\bibitem[B]{B} C. Bernardin,
%\newblock  {\it Fluctuations in the occupation time of a site in the asymmetric simple exclusion process.}
%\newblock Ann. Probab. {\bf 32} (2004), no. 1B, 855--879.
%
%\bibitem[BG]{BG} L. Bertini and G. Giacomin, 
%\newblock  {\it Stochastic Burgers and KPZ equations from particle systems},
%\newblock Comm. Math. Phys. {\bf 183} (1997), no. 3, 571--607. 

%
%\bibitem[BKS]{BKS} H. van Beijeren, R. Kutner and H. Spohn,
%\newblock  {\it Excess noise for driven diffusive systems},
%\newblock Phys. Rev. Lett.  {\bf 54} (1985), 2026--2029.
%
%\bibitem[BM]{BM} M. Bramson and T. Mountford,
%{\it Stationary blocking measures for one-dimensional nonzero mean exclusion processes},
%Ann. Probab. {\bf 30} (2002), no. 3, 1082--1130.
% 

\bibitem[BDJ]{BDJ} J. Baik, P. Deift, K. Johansson,  {\it On the distribution of the length of the longest increasing subsequence of random permutations}.  J. Amer. Math. Soc.  12  (1999),  no. 4, 1119--1178.

\bibitem[BR]{BR} J. Baik, E.M.Rains,  {\it Limiting distributions for a
polynuclear growth model with external sources}, J. Stat. Phys.  {\bf 100}
(2000), 523-542.
 
 \bibitem[BS]{BS} M. Bal\'azs, T. Sepp\"all\"ainen, {\it Order of current variance and diffusivity in the 
asymmetric simple exclusion process}, {\texttt{arxiv.org/math.PR/0608400}}
 
 \bibitem[CG]{CG} E. Cator, P. Groeneboom, {\it Second class particles and cube root asymptotics for Hammersley's process},  Ann. Probab. { \bf{34}}  (2006),  no. 4, 1273--1295.
 
\bibitem[FF]{FF}  P.A. Ferrari, L.R. Fontes,  {\it Current fluctuations for the asymmetric simple exclusion process}, Ann. Probab. {\bf 22 }(1994), no. 2, 820--832.
% 



\bibitem[FKS]{FKS} P. A. Ferrari, C. Kipnis, E. Saada,  {\it Microscopic structure of travelling waves in the asymmetric simple exclusion process},  Ann. Probab.  {\bf 19}  (1991),  no. 1, 226--244. 

\bibitem[FS]{FS}  P. L. Ferrari and H. Spohn, 
\newblock  {\it Scaling limit for the space-time covariance of the 
stationary totally asymmetric simple exclusion process}, 
%\newblock \texttt{arXiv:math-phys/0504041},
\newblock Comm. Math. Phys. {\bf 265} (2006), no. 1, 1--44

%

% 
% \bibitem[FNS]{FNS} D. Forster, D. Nelson and M. J. Stephen, 
%\newblock {\it Large-distance and long time properties of a randomly stirred fluid}, 
%\newblock Phys. Rev. A  {\bf 16} (1977), 732--749.

\bibitem[J]{J} K. Johansson, {\it Shape fluctuations and random matrices},
Comm. Math. Phys. {\bf 242} (2003), 277-329. 
 
\bibitem[KL]{KL} C. Kipnis and C. Landim, \newblock Scaling limits of interacting particle systems. Grundlheren der Mathematischen Wissenschaften 320, Springer-Verlag, Berlin, New York, 1999.
 
% 
% \bibitem[KPZ]{KPZ} K. Kardar, G. Parisi and Y. Z. Zhang, 
%\newblock {\it Dynamic scaling of growing interfaces},
%\newblock Phys. Rev. Lett.  {\bf 56} (1986), 889--892.
%
%\bibitem[L]{L}  T.M. Liggett, 
%\newblock Interacting particle systems. Grundlehren der Mathematischen Wissenschaften, 276. Springer-Verlag, New York, 1985.
%
\bibitem[LOY]{LOY}  C. Landim, S. Olla, H.T. Yau, 
\newblock {\it  Some properties of the diffusion coefficient for asymmetric simple exclusion processes.}
\newblock Ann. Probab. 24 (1996), no. 4, 1779--1808.

\bibitem[LQSY]{LQSY} C. Landim, J. Quastel, M. Salmhofer and H.-T. Yau, {\it Superdiffusivity of asymmetric exclusion process in dimensions one and two}, Comm. Math. Phys. {\bf 244} (2004), no. 3, 455--481. 

\bibitem[LY]{LY} C. Landim and H.-T. Yau, \newblock {\it Fluctuation-dissipation equation of asymmetric simple exclusion processes},\newblock  Probab. Theory Related Fields {\bf 108} (1997), no. 3, 321--356.

\bibitem[PS]{PS} M. Pr\"ahofer and H. Spohn, {\it Current fluctuations for the totally asymmetric simple exclusion process}, In and out of equilibrium (Mambucaba, 2000), 185--204, Progr. Probab. {\bf 51}, BirkhŠuser Boston, Boston, MA, 2002. 

%\bibitem[S]{S}
%S. Sethuraman
%\newblock {\it An equivalence of $H_{-1}$ norms for the 
%simple exclusion process},
%\newblock Ann. Prob. {\bf 31} (2003), No. 1, 35--62,.
% 
% 
% \bibitem[SX]{SX} S. Sethuraman and L. Xu, 
%{\it A central limit theorem for reversible exclusion and zero-range particle systems},
%Ann. Probab. {\bf 24} (1996), no. 4, 1842--1870.



\bibitem[QV]{QV} J. Quastel, B. Valk\' o, {\it $t^{1/3}$ Superdiffusivity of Finite-Range Asymmetric Exclusion Processes on $\mathbb Z$}, to appear in Comm. Math. Phys. 

\bibitem[S]{Sp} H. Spohn, Large Scale Dynamics of Interacting Particles,
Springer-Verlag, 1991.

\bibitem[V]{V} S. R. S. Varadhan,  {\it Lectures on hydrodynamic scaling}, Hydrodynamic limits and related topics (Toronto, ON, 1998), 3--40, Fields Inst. Commun. {\bf 27}, Amer. Math. Soc., Providence, RI, 2000. 
\bibitem[Y]{Y} H.-T. Yau,  {\it $(\log t)\sp {2/3}$ law of the two dimensional asymmetric simple exclusion process}, Ann. of Math. (2) {\bf 159} (2004), no. 1, 377--405. 
\end{thebibliography}
%\bibliographystyle{alpha}

\end{document}